\documentclass[reqno,dvips,12pt]{article} 

\usepackage{epsfig,graphicx,amsmath,amssymb,amsfonts}
cm \textwidth=16.0 true cm \emergencystretch=10pt

\newtheorem{theorem}{\qquad Теорема}
\newtheorem{lemma}{\qquad Lemma}

\newcommand{\cA}{{\cal A}}
\newcommand{\cM}{{\cal M}}

\newcommand{\cN}{{\cal N}}

\newcommand{\nn}{\nonumber}

\newcommand{\pa}{\partial}

\DeclareMathOperator{\re}{re}

\begin{document}

\title {\textbf{On the Minimization of Operational Risks}}

\author{\textbf{V.P.Maslov}\thanks{Moscow State University, Physics Department}}

\date{ }

\maketitle

\begin{abstract}
We give a risk-minimizing formula for government investments
taking into account the zero intelligence law for financial
markets.
\end{abstract}

A problem  of reducing  operational  risks has the specific
feature that the enterprises $A_1,A_2, \dots, A_s$ in which the
government is willing to invest are specified in advance,
possibly for political reasons, on the basis of their
profitability, or in accordance with general economic strategies.
Here $s$ is sufficiently large. Moreover, the government orders
the enterprises by priority: for each $i$, the enterprise $A_i$ is
more promising than the enterprise $A_{i-1}$. This means that one
should not buy more stock (bonds) of the enterprise $A_{i-1}$
than of $A_i$. In other respects, one should distribute the money
over these enterprises (buy stock) so as to minimize the
operational risks. Governmental structures indicate the amounts
of money to be invested in the highest priority enterprise $A_s$
and the least priority enterprise $A_1$, and also set a budget
restraint.

When solving this problem, we shall adhere to the concept of
\textit{zero intelligence}   in financial markets. This concept
arose in a study of the London Stock Exchange, and since then the
literature discussing the phenomenon has been steadily growing.
However, the zero intelligence phenomenon is quite natural from
our viewpoint \cite{DAN_Fond}.

The effect is that traders use stock exchange data to form their
portfolios at random ("without resorting to intelligence").
However, according to Kolmogorov, randomness is none other than
large complexity. For example, it does not matter if we bet on
heads or tails randomly or construct a long algorithm and stick
to it: in any case, the coin lands on heads (and, accordingly,
tails) approximately half of the times.

Let us continue the analogy. We toss a coin $S$ times, where $S$
is large, and repeat the series of $S$ tossings $N$ times. All
possible outcomes, including "all heads" or "alternating heads
and tails," are equiprobable (and have a tiny probability). But a
majority of outcomes have approximately equal numbers of heads
and tails. In particular, this means that if we play heads and
tails and bet on heads and tails half of the times each, then the
risk (and also the payoff) is minimal.

Consequently, returning to our problem, if we find the largest
"cluster" of outcomes that are equiprobable by the zero
intelligence principle, then the corresponding investment
distribution will be safest, or risk-free.

Let $C_i$ be the number of shares bought by the government
investment agency from the $i$-th enterprise (bank). By
assumption, we have the relations
\begin{equation}\label{1}
C_{i+1} \geq C_i, \ \ i=1,2, \dots, s, \ \ C_s\leq M, \ C_1\geq K,
\end{equation}
where $M$ and $K$ are given. Hence
\begin{equation}\label{2}
C_s - C_1\leq M-K.
\end{equation}
Let $p_i$ be the stock price for the ith enterprise and $N=M-k$.
Obviously,
$$
\sum^s_{i=1} C_i p_i \leq C_s \sum^s_{i=1} p_i.
$$
Hence the budget restraint $\Phi$ satisfies the relations
\begin{equation}\label{3}
C_1 \sum_{i=1}^s p_i \leq\sum^s_{i=1} C_i p_i \leq \Phi \leq C_s
\sum^s_{i=1} p_i,
\end{equation}
i.e., under this natural inequality on $\Phi$ one has
\begin{equation}\label{3a}
\sum_{i=1}^s C_i p_i\leq \Phi.
\end{equation}
We have
$$
\sum^s_{i=1} C_i p_i\equiv  \sum^s_{i=2} (C_i -C_{i-1})
\sum^s_{j=i} p_j + C_1 \sum^s_{i=1} p_i \leq \Phi,
$$
whence it follows that
\begin{equation}\label{4a}
\sum^s_{i=2} (C_i -C_{i-1})\sum^s_{j=i} p_j \leq \Phi-C_1
\sum_{i=1}^s p_i.
\end{equation}
We set
\begin{equation}\label{4}
N_i = C_i - C_{i-1}, \ \ \lambda_i = \sum^s_{j=i} p_j.
\end{equation}
Then
\begin{equation}\label{5}
\sum_{j=2}^i N_j =C_i-C_1; \ \ \lambda_i \leq\lambda_{i-1} \leq
\dots \leq \lambda_1; \ \lambda_1=\sum_{j=1}^s p_s.
\end{equation}
We obtain
\begin{equation}\label{6}
\sum^s_{i=2} N_i \lambda_i \leq \Phi-K \sum_{j=1}^s p_s.
\end{equation}
Let $E=\Phi-\sum_{i=1}^s p_i K$.

Consider the set of all sequences  $\{N_i\}$ such that
\begin{equation}\label{7}
\sum_{i=2}^s N_i \lambda_i \leq E.
\end{equation}

Since all sequences $\{C_i\}$ satisfying conditions ~(\ref{1})
and ~(\ref{2}) are equiprobable, it follows that the sequences
$N_i=C_i-C_{i-1}; \  i=2, \dots, s,$ satisfying the relation
\begin{equation}\label{7a}
\sum_{i=2}^s N_i = M-K =N,
\end{equation}
are equiprobable under the condition $N_i \geq 0$ and condition
~(\ref{7}).

We assume that $a_1N \leq s \leq a_2N$, $a_1, a_2$ are constants.
We introduce the notation: $\cM$ is the set of all sets $\{N_i\}$
satisfying conditions~(\ref{7}) and~(\ref{7a}); \ $\cN\{\cM\}$ is
the number of elements of the set~$\cM$.

\begin{lemma} \label{theor1} 
Suppose that all the variants of sets $\{N_i\}$ satisfying the
conditions ~(\ref{7}) and ~(\ref{7a}) are equiprobable. Then the
number of variants $\cN$ of sets $\{N_i\}$ satisfying
conditions~(\ref{7}) and~(\ref{7a}) and the additional relation
\begin{equation} \label{theorema1}
|\sum^l_{i=1} N_i - \sum^l_1\frac{q_i}{e^{\beta'
\lambda_i-\nu'}-1}|\geq  N^{(3/4+\varepsilon)}
\end{equation}
is less than $\frac{c_1 \cN\{\cM\}}{N^m}$ (where~$c_1$ and~$m$
are any arbitrary numbers, $\sum_{i=1}^l q_i \geq\varepsilon Q$,
and $\varepsilon$ is arbitrarily small).
\end{lemma}

{\bf{\qquad Proof of Lemma 1.}}

Let $\cA$ be a subset of $\cM$ satisfying the condition
$$
|\sum_{i=l+1}^s N_i - \sum_{i=l+1}^s \frac{q_i}
{e^{\beta\lambda_i-\nu}-1}|\leq \Delta;
$$
$$
|\sum_{i=1}^l N_i-\sum_{i=1}^l \frac{q_i}
{e^{\beta'\lambda_i-\nu'}-1}|\leq \Delta,
$$
where $\Delta$, $\beta$, $\nu$ are some real numbers independent
of~$l$.

We denote
$$
|\sum_{i=l+1}^s N_i-\sum_{i=l+1}^s \frac{q_i}
{e^{\beta\lambda_i-\nu}-1}| =S_{s-l};
$$
$$
|\sum_{i=1}^l N_i-\sum_{i=1}^l \frac{q_i}
{e^{\beta'\lambda_i-\nu'}-1}| =S_l.
$$

Obviously, if $\{N_i\}$ is the set of all sets of integers on the
whole, then
\begin{equation}\label{Proof1}
\cN\{\cM \setminus \cA\} = \sum_{\{N_i\}} \Bigl(
\Theta(E-\sum_{i=1}^s N_i\lambda_i) \delta_{(\sum_{i=1}^s N_i),N}
\Theta(S_l-\Delta)\Theta (S_{s-l}-\Delta)\Bigr),
\end{equation}
where $\sum N_i=N$.

Here the sum is taken over all integers $N_i$, $\Theta(\lambda)$
is the Heaviside function, and $\delta_{k_1,k_2}$ is the
Kronecker symbol.

We use the integral representations
\begin{eqnarray}
&&\delta_{NN'}=\frac{e^{-\nu N}}{2\pi}\int_{-\pi}^\pi d\varphi
e^{-iN\varphi} e^{\nu N'}e^{i N'\varphi},\label{D7}\\
&&\Theta(y)=\frac1{2\pi i}\int_{-\infty}^\infty
d\lambda\frac1{\lambda-i}e^{\beta y(1+i\lambda)}.\label{D8}
\end{eqnarray}

Now we perform the standard regularization. We replace the first
Heaviside function~$\Theta$ in~(\ref{Proof1}) by the continuous
function
\begin{equation}
\Theta_{\alpha}(y) =\left\{
\begin{array}{ccc}
0 &\mbox{for}& \alpha > 1, \  y<0 \nn \\
1-e^{\beta y(1-\alpha)} &\mbox{for}& \alpha > 1,\  y \geq 0,
\label{Naz1}
\end{array}\right.
\end{equation}
\begin{equation}
\Theta_{\alpha}(y) =\left\{
\begin{array}{ccc}
e^{\beta y(1-\alpha)} &\mbox{for}&\alpha < 0, \ y<0 \nn \\
1 &\mbox{for}& \alpha < 0, \ y \geq 0, \label{Naz2}
\end{array}\right.
\end{equation}
where $\alpha \in (-\infty,0) \cup (1, \infty)$ is a parameter,
and obtain
\begin{equation}\label{proof2}
\Theta_\alpha(y) = \frac1{2\pi i} \int_{-\infty}^{\infty}
e^{\beta y(1+ix)} (\frac 1{x-i} - \frac 1{x-\alpha i}) dx.
\end{equation}

If  $\alpha > 1$, then $\Theta(y)\leq \Theta_{\alpha}(y)$.

Let $\nu <0$. We substitute~(\ref{D7}) and~(\ref{D8})
into~(\ref{Proof1}),  interchange the integration and summation,
then pass to the limit as $\alpha \to \infty$ and obtain the
estimate
\begin{eqnarray}
&&\cN\{\cM \setminus   \cA\} \leq \nn \\
&&\leq \Bigl|\frac{e^{-\nu N+\beta E }}{i(2\pi)^2}\int_{-\pi}^\pi
\bigl[ \exp(-iN\varphi)
\sum_{\{N_j\}}\Bigl(\exp\bigl\{(-\beta\sum_{j=1}^s
N_j\lambda_j)+(i\varphi+\nu)
N_j\bigr\}\bigr]\ d\varphi \times \nn \\
&& \times \Theta(S_l -\Delta)\Theta(S_{s-l}-\Delta)\Bigr|, \ \sum
N_i=N,
\end{eqnarray}
where $\beta$ and $\nu$ are real parameters such that the series
converges for them.

To estimate the expression in the right-hand side, we bring the
absolute value sign inside the integral sign and then inside the
sum sign, integrate over $\varphi$, and obtain
\begin{eqnarray}
&&\cN\{\cM \setminus \cA\} \leq \frac{e^{-\nu N+\beta E }}{2\pi}
\sum_{\{N_i\}}\exp\{-\beta\sum_{i=1}^sN_i\lambda_i+\nu
N_i\}\times \nn \\
&& \times\Theta (S_l-\Delta)\Theta (S_{s-l}-\Delta).
\end{eqnarray}

We denote
\begin{equation}\label{D9a}
Z(\beta,N)=\sum_{\{N_i\}} e^{-\beta\sum_{i=1}^s N_i\lambda_i},
\end{equation}
where the sum is taken over all $N_i$ such that $\sum_{i=1}^s
N_i=N$,
$$
\zeta_l(\nu,\beta)= \prod_{i=1}^{l} \xi_i\left(\nu,\beta\right);
\zeta_{s-l}(\nu,\beta)= \prod_{i=l+1}^{s}
\xi_i\left(\nu,\beta\right);
$$
$$
\quad \xi_i(\nu,\beta)=
\frac{1}{(1-e^{\nu-\beta\lambda_i})^{q_i}}, \qquad i=1,\dots,l.
$$

It follows from the inequality for the hyperbolic cosine
$\cosh(x)=(e^x+e^{-x})/2$ for $|x_1| \geq \delta; |x_2| \geq
\delta$:
\begin{equation}
\cosh(x_1)\cosh(x_2)> \frac{e^\delta}{2} \label{D33}
\end{equation}
that the inequality
\begin{equation}
\Theta(S_{s-l}-\Delta) \Theta(S_{l}-\Delta)\le e^{-c\Delta}
\cosh\Bigl(c\sum_{i=1}^{l} N_i
-c\phi_l\Bigr)\cosh\Bigl(c\sum_{i=l+1}^{s} N_i
-c\overline{\phi}_{s-l}\Bigr), \label{D34}
\end{equation}
where
$$
\phi_l= \sum_{i=1}^l \frac{q_i}{e^{\beta'\lambda_i-\nu'}-1};
\qquad \overline{\phi}_{s-l}= \sum_{i=l+1}^s
\frac{q_i}{e^{\beta\lambda_i-\nu}-1},
$$
holds for all positive $c$ and~$\Delta$.

We obtain
\begin{eqnarray}
&&\cN\{\cM \setminus \cA\} \leq  e^{-c\Delta} \exp\left(\beta
E-\nu N\right) \times \nn \\
&& \times \sum_{\{N_i\}}\exp\{-\beta\sum_{i=1}^l
N_i\lambda_i+\nu\sum_{i=1}^l N_i\} \cosh\left(\sum_{i=1}^{l} c N_i -
c\phi\right) \times \nn \\
&& \times \exp\{-\beta\sum_{i=l+1}^s N_i\lambda_i +\nu
\sum_{i=l+1}^s
N_i\} \cosh\Bigl(\sum_{i=l+1}^s c N_i -c\overline{\phi}\Bigr) = \nn \\
&& =e^{\beta E} e^{-c\Delta} \times \nn \\
&&\times \left( \zeta_l(\nu-c,\beta) \exp(-c\phi_{l})
+\zeta_l(\nu+c,\beta)\exp(c\phi_{l})\right) \times \nn \\
&&\times\left(\zeta_{s-l}(\nu-c,\beta)
\exp(-c\overline{\phi}_{s-l})
+\zeta_{s-l}(\nu+c,\beta)\exp(c\overline{\phi}_{s-l})\right).
\label{5th}
\end{eqnarray}

Now we use the relations
\begin{equation}\label{5tha}
\frac {\pa}{\pa\nu}\ln \zeta_l|_{\beta=\beta',\nu=\nu'}\equiv
\phi_l; \frac {\pa}{\pa\nu}\ln
\zeta_{s-l}|_{\beta=\beta',\nu=\nu'}\equiv \overline{\phi}_{s-l}
\end{equation}
and the expansion $\zeta_l(\nu\pm c,\beta)$ by the Taylor formula.
There exists a $ \gamma <1$ such that
$$
\ln(\zeta_l(\nu\pm c,\beta)) =\ln\zeta_l(\nu,\beta)\pm
c(\ln\zeta_l)'_\nu(\nu,\beta)+\frac{c^2}{2}(\ln\zeta_l)^{''}_\nu
(\nu\pm\gamma c,\beta).
$$
We substitute this expansion, use formula~(\ref{5tha}), and see
that $\phi_{\nu,\beta}$ is cancelled.

Another representation of the Taylor formula implies
\begin{eqnarray}
&&\ln\left(\zeta_l(\nu+c,\beta)\right)=
\ln\left(\zeta_l(\beta,\nu)\right)+
\frac{c}\beta\frac{\pa}{\pa\nu}\ln\left(\zeta_l(\beta,\nu)\right)+\nn\\
&&+\int_{\nu}^{\nu+c/\beta}d\nu' (\nu+c/\beta-\nu')
\frac{\pa^2}{\pa\nu'^2}\ln\left(\zeta_l(\beta,\nu')\right).\label{CC1}
\end{eqnarray}
A similar expression holds for $\zeta_{s-l}$.

From the explicit form of the function $\zeta_l(\beta,\nu)$, we
obtain
\begin{equation}
\frac{\pa^2}{\pa\nu^2}\ln\left(\zeta_l(\beta,\nu)\right)=
\beta^2\sum_{i=1}^{l}
\frac{g_i\exp(-\beta(\lambda_i+\nu))}{(\exp(-\beta(\lambda_i+\nu))-1)^2}\leq
\beta^2Qd, \label{CC2}
\end{equation}
where $d$ is given by the formula
$$
d=\frac{\exp(-\beta(\lambda_1+\nu))}{(\exp(-\beta(\lambda_1+\nu))-1)^2}..
$$
The same estimate holds for $\zeta_{s-l}$.

Taking into account the fact that $\zeta_l\zeta_{s-l}=\zeta_s$,
we obtain the following estimate for $\beta=\beta'$ and
$\nu=\nu'$:
\begin{equation}\label{eval1}
\cN\{\cM \setminus \cA\}
\leq\zeta_s(\beta',\nu')\exp(-c\Delta+\frac{c^2}{2}\beta^2Qd)
\exp(E\beta'-\nu'N).
\end{equation}

Now we express $\zeta_s(\nu',\beta')$ in terms $Z(\beta,N)$. To
do this, we prove the following lemma 2.

\begin{lemma}
Under the above assumptions, the asymptotics of the integral
\begin{equation}\label{lemma_1}
Z(\beta,N) = \frac{e^{-\nu N}}{2\pi}\int_{-\pi}^\pi
 d\alpha e^{-iN\alpha}\zeta_s(\beta,\nu+i\alpha)
\end{equation}
has the form
\begin{equation}\label{lemma_2}
Z(\beta,N) = C e^{-\nu N} \frac{\zeta_s(\beta,\nu)}{|(\partial^2
\ln\zeta_s(\beta,\nu))/ (\partial^2\nu)|} (1+O(\frac 1N)),
\end{equation}
where $C$ is a constant.
\end{lemma}

We have
\begin{equation}
 Z(\beta,N) = \frac{e^{-\nu N}}{2\pi}\int_{-\pi}^\pi
 e^{-iN\alpha}\zeta_s(\beta,\nu+i\alpha)\,d\alpha
 =\frac{e^{-\nu N}}{2\pi}\int_{-\pi}^\pi e^{NS(\alpha,N)} d\alpha ,\label{D15}
\end{equation}
where
\begin{equation}\label{qq}
    S(\alpha,N) = -i\alpha+ \ln \zeta_s (\beta, \nu +i\alpha)
    = -i\alpha - \sum_{i=1}^s q_i\ln [1-e^{\nu+i\alpha-\beta\lambda_i}].
\end{equation}
Here $S$ depends on $N$, because $s$, $\lambda_i$, and $\nu$ also
depend on~$N$; the latter is chosen so that the point $\alpha=0$
be a stationary point of the phase~$S$, i.e., from the condition
\begin{equation}\label{qq1}
 N=\sum_{i=1}^s\frac{q_i}{e^{\beta\lambda_i-\nu}-1}.
\end{equation}
We assume that $a_1N \leq s \leq a_2N$, $a_1,
a_2=\operatorname{const}$, and, in addition, $0\le\lambda_i\le B$
and $B=\operatorname{const}$, $i=1,\dots,s$.
 If these conditions are satisfied
in some interval $\beta\in[0,\beta_0]$ of the values of the
inverse temperature, then all the derivatives of the phase are
bounded, the stationary point is nondegenerate, and the real part
of the phase outside a neighborhood of zero is strictly less than
its value at zero minus some positive number. Therefore,
calculating the asymptotics of the integral, we can replace the
interval of integration $[-\pi,\pi]$ by the interval
$[-\varepsilon,\varepsilon]$. In this integral, we perform the
change of variable
\begin{equation}\label{qqq}
  z=\sqrt {S(0,N)-S(\alpha,N)}.
\end{equation}
This function is holomorphic in the disk
$|{\alpha}|\le\varepsilon$ in the complex $\alpha$-plane and has
a holomorphic inverse for a sufficiently small~$\varepsilon$. As
a result, we obtain
\begin{equation}\label{qqq1}
    \int_{-\varepsilon}^\varepsilon e^{NS(\alpha,N)}
    d\alpha=e^{NS(0,N)}\int_\gamma e^{-Nz^2}f(z)\,dz,
\end{equation}
where the path $\gamma$ in the complex $z$-plane is obtained from
the interval $[-\varepsilon, \varepsilon]$ by the
change~\eqref{qqq} and
\begin{equation}\label{qqq3}
    f(z)=\left(\frac{\partial\sqrt {S(0,N)-S(\alpha,N)}}
    {\partial\alpha}\right)^{-1}\bigg|_{\alpha=\alpha(z)}.
\end{equation}
For a small~$\varepsilon$ the path $\gamma$ lies completely
inside the double sector $\re(z^2)>c(\re z)^2$ for some $c>0$;
hence it can be ``shifted'' to the real axis so that the integral
does not change up to terms that are exponentially small in~$N$.
Thus, with the above accuracy, we have
\begin{equation}\label{qqq4}
  Z(\beta,N) =  \frac{e^{-\nu N}}{2\pi}\int_{-\varepsilon}^\varepsilon e^{-Nz^2}f(z)\,dz.
\end{equation}
Since the variable $z$ is now real, we can assume that the
function $f(z)$ is finite (changing it outside the interval of
integration), extend the integral to the entire axis (which again
gives an exponentially small error), and then calculate the
asymptotic expansion of the integral expanding the integrand in
the Taylor series in~$z$ with a remainder. This justifies that
the saddle-point method can be applied to the above integral in
our case. The proof of the lemma 1 is complete.

Now we estimate $Z(\beta,N)$. To do this, we proof  the following
lemma 3.

\begin{lemma}
The quantity
\begin{equation}\label{qqq5}
\frac{1}{\cN(\cM)} \sum_{\{N_i\}} e^{-\beta\sum_{i=1}^s
N_i\lambda_i},
\end{equation}
where $\sum N_i =N$ and $\lambda_iN_i\leq E-N^{1/2+\varepsilon}$,
tends to zero faster than $N^{-k}$ for any $k$, $\varepsilon>0$.
\end{lemma}

We consider the point of minimum in $\beta$ of the right-hand
side of ~(\ref{5th}) with $\nu(\beta,N)$ satisfying the condition
$$
\sum \frac{q_i}{e^{\beta\lambda_i-\nu(\beta,N)}-1} =N.
$$
It is easy to see that it satisfies condition~(\ref{7}). Now we
assume that the assumption of the lemma is not satisfied.

Then for $\sum N_i=N$,  $\sum \lambda_i N_i\geq
E-N^{1/2+\varepsilon}$, we have
$$
e^{\beta E}\sum_{\{N_i\}} e^{-\beta\sum_{i=1}^s N_i\lambda_i}\geq
e^{(N^{1/2}+\varepsilon)\beta}.
$$
Obviously, $\beta\ll \frac{1}{\sqrt{N}}$ provides a minimum
of~(\ref{5th}) if the assumptions of Lemma~1 are satisfied, which
contradicts the assumption that the minimum in~$\beta$ of the
right-hand side of~(\ref{5th}) is equal to~$\beta'$.

We set $c=\frac\Delta{N^{1+\alpha}}$ in formula~(\ref{eval1})
after the substitution~(\ref{lemma_2}); then it is easy to see
that the ratio
$$
\frac{\cN(\cM \setminus\cA)}{\cN(\cM)}\approx \frac 1{N^m},
$$
where $m$ is an arbitrary integer, holds for
$\Delta=N^{3/4+\varepsilon}$. The proof of the lemma 1 is
complete.

From ~(\ref{4}) - ~(\ref{6}) follows

\begin{theorem}
A risk-free investment (stock purchase) under conditions
~(\ref{1})-- (\ref{3}) takes place if
\begin{equation}\label{10b}
C_i= C_1 + \sum_{j=2}^i \frac{1}{e^{\beta \sum_{k=j}^s
p_k-\sigma}-1}+O(N^{3/4}); \ i>\delta N, \ \delta>0,
\end{equation}
where $\sigma$ and $\beta$  are determined from the conditions
$$
\sum_{j=2}^s \frac{1}{e^{\beta \sum_{k=j}^s p_k-\sigma}-1}=M-K; \
\ \sum_{i=2}^s p_i\sum_{j=2}^i \frac{1}{e^{\beta\sum_{k=j}^s p_k
-\sigma}-1} =\Phi - p_1K.
$$
\end{theorem}


\begin{thebibliography}{99}
\bibitem{DAN_Fond}
V. ~P. ~Maslov, ''On the increasing complexity principle in
portfolio formation at stock exchange'', Dokl. Ross. Akad. Nauk
[Russian Acad. Sci. Dokl. Math.], 404 (2005), no. 4, 446-450.
\bibitem{NelinSred}
V.~P.~Maslov. The nonlinear average in economics.
// Mat. Zametki [Math. Notes]. 2005, \textbf{78}, No.~3, 377--395.


\end{thebibliography}
\end{document}